\documentclass[12pt]{article}

\usepackage[english]{babel}
\usepackage{graphicx}
\usepackage{latexsym}
\usepackage{amsfonts}
\usepackage{amscd}
\usepackage{amssymb}
\usepackage{amsmath}
\usepackage{mathrsfs}
\usepackage{amsthm}

\input xy
\xyoption{all}

\DeclareMathOperator{\supp}{supp}

\newcommand{\re}{\mathbb{R}}
\newcommand{\s}{\mathcal S}
\newcommand{\G}{\mathcal G}

\newcommand{\rn}{\mathbb{R}^n}
\newcommand{\rp}{\mathcal{R}_{punc}}

\newcommand{\op}{\Omega_{punc}}

\newtheorem{lem}{Lemma}
\newtheorem{theo}{Theorem}
\newtheorem{defn}{Definition}

\newtheorem{prop}{Proposition}

\def\Wagner{W}
\def\Preston{P}
\def\HZ{HZ}
\def\Barnes{B}
\def\Munn{M}
\def\DP{DP}
\def\Renault{R}
\def\Exel{E}

\begin{document}
\title{The tiling C*-algebra viewed as a tight inverse semigroup algebra}
\author{Ruy Exel\thanks{Departamento de Matem\'atica, Universidade
Federal de Santa Catarina, 88040-900 Florian\'opolis SC, Brazil. Email: \texttt{exel@mtm.ufsc.br}. Partially supported by CNPq.}, \and Daniel Gon\c{c}alves\thanks{ Departamento de Matem\'atica, Universidade
Federal de Santa Catarina, 88040-900 Florian\'opolis SC, Brazil. E-mail: \texttt{daemig@gmail.com}.}, \and Charles Starling\thanks{Department of Mathematics, University of Ottawa, 585 King Edward, Ottawa ON, K1N 6N5. Email: \texttt{cstar050@uottawa.ca}.}}
\date{}
\maketitle
\begin{abstract}
We realize Kellendonk\'{}s C*-algebra of an aperiodic tiling as the tight C*-algebra of the inverse semigroup associated to the tiling, thus providing further evidence that the tight C*-algebra is a good candidate to be the natural associative algebra to go along with an inverse semigroup.
\end{abstract}

\section{Introduction}

The notion of an inverse semigroup was introduced independently by Wagner
\cite{\Wagner} and Preston \cite{\Preston} in the early 50's to give an abstract
formulation for the algebraic structure formed by all partialy defined bijective maps on a given
set, under the operation of composition in the largest domain where it makes
sense.
  Perhaps the best way to illustrate that inverse semigroups do yield such an
abstraction is via what is now known as the Wagner--Preston Theorem, according
to which every inverse semigroup is isomorphic to a semigroup of partially
defined bijections.

One of the main themes behind the present work is the question of what is the
natural associative algebra to go along with a given inverse semigroup, generalizing
the notion of the group algebra of a given group.
Since inverse semigroups are special examples of semigroups, they may of course
be studied simply as semigroups.  In particular 
  one may construct their semigroup algebra, as first studied by Munn \cite{\Munn},
  or their $l_1$-algebra as discussed by Hewitt and Zuckerman \cite{\HZ}.
  The $l_1$-algebra of an inverse semigroup has been considered by Barnes
\cite{\Barnes}, where the main problem treated in \cite{\HZ} was solved in a
very natural way.

The $l_1$-algebra of an inverse semigroup over the field of the complex numbers
has a natural structure of a *-algebra and hence it is natural to consider its
envelopping C*-algebra, as studied by Duncan and Paterson \cite{\DP}.

The extensive literature on inverse semigroup algebras  (be it the purely
algebraic object, the $l_1$-algebra, or  the envelopping C*-algebra) could be
considered definitive were it not for the fact that it has some limitations.
The first hint that the established notion fails to deliver in some situatons  has
already been noticed more than thirty years ago in Renault's thesis
\cite[2.8]{\Renault}.  
  For example, if one lets ${\mathscr{ O}}_n$ be the ``Cuntz inverse semigroup''
\cite[2.2]{\Renault}, it is perhaps not unreasonable to expect that the
C*-algebra of ${\mathscr{ O}}_n$ be isomorphic to the Cuntz algebra ${\mathcal O}_n$.
Unfortunately, this is not so: Duncan and Paterson's construction applied to
${\mathscr{ O}}_n$ gives the Toeplitz extension of ${\mathcal O}_n$, rather than ${\mathcal
O}_n$ itself.

The trouble is that the famous Cuntz relation ``$\sum_{i=1}^nS_iS_i^* = 1$"
involves sums,  besides products,  and hence it is not obvious how to express it
in (multiplicative) semigroup terms.

Addressing this question the first named author has introduced an alternative
construction of a C*-algebra from an inverse semigroup, based on the idea of
``tight representations'', which, when applied to the Cuntz inverse semigroup,
does produce the Cuntz algebra.

The purpose of the present work is to apply these ideas to inverse semigroups
which naturally occur in the study of tilings, the study of which traces back to Wang \cite{Wang}. A tiling of $\re^n$ gives rise to an inverse semigroup, which was originally introduced and studied by Kellendonk \cite{K, Kel3}: Given a tiling $T$ of $\re^n$, one considers the set of all triples
$(t_1,P,t_2)$, where $P$ is a patch in $T$ (see below for the precise
definitions), and $t_1$ and $t_2$ are tiles in $P$.  The equivalence class of
such a triple modulo the action of the group of translations of $\re^n$ is often
called a ``doubly-pointed pattern class''.
  The set $\mathcal S$ of all doubly-pointed pattern classes, to which is added an ad-hoc
``zero element'', forms an inverse semigroup under the multiplication defined by
  $$
  [t_1,P,t_2] [t_1',P',t_2'] =   [t_1,P \cup P',t_2'],
  $$
  provided $t_2=t_1$ and $P\cup P'$ is a patch in $T$. When these conditions
are not met the product is simply set to be zero.

In a close parallel to the above question of whether or not one obtains the
Cuntz algebra as the C*-algebra of the Cuntz inverse semigroup, one may ask if
the C*-algebra of an aperiodic tiling $T$, as defined in \cite{Kel2}, may be
obtained from the inverse semigroup $\mathcal S$ described above. The C*-algebra associated to a tiling is the C*-algebra of an \' etale equivalence relation $\rp$ which is translational equivalence on a suitable space of tilings. Under suitable conditions on the tilings considered, the algebra $C^*(\rp)$ contains a generating set which form an inverse semigroup with respect to the C*-algebra product which is isomorphic to the inverse semigroup associated to a tiling above. 

The answer to whether $C^*(\rp)$ can be obtained from $\mathcal S$ is again negative for Duncan and Paterson's construction. Kellendonk addresses this issue for the specific example of a tiling in \cite{Kel3} and for general 0-$E$-unitary inverse semigroups in Lawson's book \cite[9.2]{Lawson}. As our main result
shows,  the answer also becomes affirmative if we use the ``tight C*-algebra'' of $\mathcal S$ instead. 

This may be viewed as further evidence that the tight C*-algebra of an inverse
semigroup is a good candidate to be the natural associative
algebra to go along with an inverse semigroup.

\section{Inverse Semgroups}

Recall that a semigroup $\mathcal{S}$ is said to be an {\bf inverse semigroup} if for every $s\in\mathcal{S}$ there is an element $s^*\in\mathcal{S}$ such that
\[
s^*ss^* = s^* \hspace{0.5cm}\textnormal{and}\hspace{0.5cm} ss^*s = s.
\]
We denote by $E$ the set of idempotent elements in $\mathcal{S}$. 

Let $X$ be a set. Denote by $\mathcal{I}(X)$ the inverse semigroup of all bijections between subsets of $X$, with the operation of composition of functions in the largest domain on which the composition makes sense. Let $\mathcal{S}$ be an inverse semigroup and let $X$ be a locally compact Hausdorff space. An {\bf action} of $\mathcal{S}$ on $X$ is a semigroup homomorphism
\[
\theta: \mathcal{S} \to \mathcal{I}(X)
\]
such that for every $s\in\mathcal{S}$ the map $\theta_s$ is continuous and its domain is open in $X$, and the union of the domains of the $\theta_s$ is all of $X$. The object we use to study actions of inverse semigroups is the {\bf groupoid of germs} associated the action. 

Let $\mathcal S$ be an inverse semigroup acting on a locally compact Hausdorff space $X$ by an action $\theta$. We denote by $Y$ the subset of $\s \times X$  given by $Y = \{ (s,x) \in \s \times X: x \in D_{s^*s} \}$, where $D_{s^*s}$ denotes the domain of $\theta_{s^*s}$, and for every $(s,x)$ and $(t,y)$ in $Y$ we say that $(s,x) \sim (t,y)$, if $x=y$ and there exists an idempotent $e$ such that $x \in D_e$ and $se=te$. We denote the equivalence class of $(s,x)$ by $[s,x]$ and call it the germ of $s$ at $x$. The groupoid of germs is then defined as $\G= Y /\sim$, with the set of composible pairs given by $\G^{\left(2 \right)}=\{ \left( [s,x],[t,y]\right) \in \G \times \G:x = \theta_t(y) \}$ and operations defined by $[s,x] \cdot [t,y] = [st,y]$, $[s,x]^{-1}=[s^*,\theta_s(x)]$. A basis for the topology of $\G$ consists or the collection of all $\Theta(s,U)$, where $\Theta(s,U)=\{[s,x]\in \G: x \in U\}$, for any $s \in \s$ and any open subset $U\subseteq D_{s^*s}$.

\section{Tilings}

A {\bf tile} is a subset of $\rn$ homeomorphic to the closed unit ball. A {\bf partial tiling} is a collection of tiles whose interiors are pairwise disjoint. A finite partial tiling will be called a {\bf patch}. The {\bf support} of a partial tiling is the union of its tiles. We define a {\bf tiling} to be a partial tiling whose support is $\rn$. Given $U\subset\rn$ and a partial tiling $T$, $T(U)$ is all the tiles that intersect $U$, that is, $T(U) = \{t\in T \mid t\cap U\neq \emptyset\}$. For $x\in\rn$, $T(\{x\})$ is frequently abbreviated as $T(x)$. Two partial tilings $T$ and $T'$ are said to {\bf agree on $U$} if $T(U) = T'(U)$, and are said to be {\bf compatible} if $T$ and $T'$ agree on Int$(\supp(T)\cap\supp(T'))$.

Given a vector $x\in\rn$ we can take any subset $U\subset \rn$ and form its translate by $x$, namely $U+x = \{ u+x \mid u\in U\}$. Thus, given a tiling $T$ we can form another tiling by translating every tile by $x$. We denote the new tiling by
\[
T+x = \{t+x \mid t\in T\}
\]

A set of tiles $\mathcal{P} = \{p_1, p_2, \dots, p_N\}$ is called a set of {\bf prototiles} for $T$ if $t\in T \Rightarrow t = p_i + x$ for some $p_i\in\mathcal{P}$ and $x\in\rn$. Prototiles may carry labels to distinguish translationally equivalent tiles. We also insist that the $p_i$, when viewed as subsets of $\rn$, have the origin in their interior. This allows us to define a designated point in each tile. If $t$ is a tile in a tiling and $t = p + x(t)$ for some $p\in\mathcal{P}$ we say that the {\bf puncture} of $t$ is $x(t)$.  

A tiling for which $T+x = T$ for some non-zero $x\in\rn$ is called {\bf periodic}. A tiling for which no such non-zero vector exists is called {\bf aperiodic}.

Given a set $\mathcal{T}$ of tilings, we would define a metric on $\mathcal{T}$. This definition is standard, for instance see \cite{KP}.

\begin{defn}
Suppose that $\mathcal{T}$ is a set of tilings and that $T, T'\in \mathcal{T}$. We define the distance between $T$ and $T'$ to be
\begin{eqnarray*}
d(T,T') & = &  \inf\{1,\epsilon \mid \exists\ x, x'\in\re^n \ni \left|x\right|, \left|x'\right| < \epsilon,\\
        &   & \  \ (T - x)(B(0,1/\epsilon)) = (T'-x')(B(0,1/\epsilon))\}.
\end{eqnarray*} 
This is called the {\bf tiling metric}.
\end{defn}
Two tilings are close in this metric if they agree on a large ball around the origin up to a small translation. If $T$ is a tiling, then we complete the set $T+\rn = \{T+x \mid x\in\rn\}$ in this norm and obtain a metric space $\Omega_T$. This space is called the {\bf continuous hull} of $T$ (also called the {\bf tiling space} associated with $T$). One can verify that given a Cauchy sequence $\{T_i\}_{i\in\mathbb{N}}$ in $T+\rn$ one can find a tiling $T'$ with $T_n \to T'$, and that if $\mathcal{P}$ is a set of prototiles for $T$ then $\mathcal{P}$ is a set of prototiles for $T'$. The space $\Omega_T$ can also be characterized as the set of all tilings $T'$ such that every patch in $T'$ appears as some translate of a patch in $T$.

Although the inverse semigroup of a tiling defined below can be defined in a more general setting, the C*-algebra of a tiling has most interest under some additional assumptions on the tiling.
\begin{defn}{\em (Assumption 1)}
A tiling $T$ is said to have {\bf finite local complexity} if for every $r>0$ the set $\{ T(B(x,r)) \mid x\in\rn \}$ contains only finitely many different patches modulo translation.
\end{defn}
If $T$ has finite local complexity, then $\Omega_T$ is compact; see \cite{RW}, Lemma 2.
\begin{defn}{\em (Assumption 2)}
A tiling $T$ is said to have {\bf repetitivity} (or is {\bf repetitive}) for any patch $P\subset T$ there exists $R>0$ such that for every $x\in\rn$ there exists $y\in\re^n$ such that $P+y\subset T(B_R(x))$. 
\end{defn}
The tiling $T$ is repetitive if and only if for every $T'\in \Omega_T$ we have $\Omega_{T'} = \Omega_{T}$; see \cite{So}, Lemma 1.2.
\begin{defn}{\em (Assumption 3)}
A tiling $T$ is called {\bf strongly aperiodic} if $\Omega_T$ contains no periodic tilings.
\end{defn}
Suppose that $\mathcal{P}$ is a set of prototiles for $T$. We let $\op \subset \Omega_T$ be the set of tilings in $\Omega_T$ that have a puncture at the origin. In other words, $T\in\op$ if and only if $T(0)\in \mathcal P$. The space $\op$ is sometimes called the {\bf punctured hull} of the tiling.
If $T$ satisifies the three assumptions above, then its punctured hull is homeomorphic to a Cantor set (see for example \cite{KP}, p. 187).  The space $\op$ has a neighbourhood base consisting of sets of the following form: for a patch $P$ and tile $t\in P$, define
\[
U(P,t) = \{ T\in\op \mid P-x_t \subset T \}.
\]
Each of these sets is clopen and the set of all such $U(P,t)$ forms a neighbourhood base for $\op$. Notice that for $x\in\rn$, we have $U(P,t) = U(P+x, t+x)$. 

We now can define the inverse semigroup associated to a tiling $T$. We let $\mathcal M$ be the set of patches $P$ such that the $P$ appears somewhere in $T$ and the support of $P$ has connected interior. Let $P, Q\in\mathcal{M}$, $t_1, t_2\in P$ and $s_1, s_2\in Q$. We say the triples $(t_1, P, t_2)$ and $(s_1, Q, s_2)$ are equivalent, and write $(t_1, P, t_2)\sim (s_1, Q, s_2)$, if there exists $x\in\rn$ such that $t_1 + x = s_1$, $t_2 + x = s_2$ and $P+x = Q$. The equivalence class of $(t_1, P, t_2)$ is denoted $[t_1, P, t_2]$ and is called a ``doubly pointed pattern class''. Let
\[
\mathcal S = \{ [t_1,P, t_2] \mid P\in \mathcal {M},\hspace{0.2cm} t_1, t_2\in P\}\cup\{0\}.
\]
We give $\mathcal S$ the structure of an inverse semigroup as follows. Let $[t_1,P, t_2]$ and $[t_1',P', t_2']$ be two doubly pointed pattern classes. Suppose we can find $x, x'\in\rn$ such that $P+x$ and $P'+x'$ are patches in $T$ and such that $t_2 +x = t_1' + x'$. Then we define the product of these two classes to be
\[
[t_1, P, t_2][t_1', P', t_2'] = [t_1+x, (P+x)\cup(P'+x'), t_2'+x']
\]
This product is well-defined because $(P+x)\cup(P'+x')$ is a patch in $T$. We notice that the product $[t_1, P, t_2][t_1', P', t_2']$ is non-zero if and only if $U(P,t_2)\cap U(P', t_1') \neq \emptyset$. If no such translations exist, we define the product to be zero. We also define all products involving 0 to be 0, and we define
\[
[t_1,P, t_2]^* = [t_2, P, t_1].
\]
This product and involution give $\mathcal S$ the structure of an inverse semigroup; see for example \cite[9.5]{Lawson}, \cite{Kel3} or \cite{K}. Notice when forming the product $[t_1, P, t_2][t_1', P', t_2']$ we may pick representatives in each class so that $t_2 = t_1'$ and $P\cup P'$ is a patch. In this case we obtain the nicer formula $$[t_1, P, t_2][t_1', P', t_2'] = [t_1, P\cup P', t_2'].$$


We notice that there is an action $\theta^{\Omega}$ of $\mathcal S$ on the space $\op$. Let $s = [t_1, P, t_2]$ and define $x_s$ as $x(t_1)-x(t_2)$. Then $s$ can be seen as a partial bijection
\[
\theta^\Omega_s: U(P, t_2) \to U(P, t_1)
\]
\[
\theta^\Omega_s(T) = T -x_s.
\]

When we restrict to $\op$, we no longer have an $\rn$ action. If $T$ is a punctured tiling, then $T+x$ need not be a a punctured tiling. We may still consider translational equivalence restricted to this subspace.

Define
\[
\rp = \{ (T, T+x) \mid T, T+x \in \op \}
\]
If we give this the topology from $\op \times \rn$, then $\rp$ becomes a principal $r$-discrete groupoid. The unit space of $\rp$ is clearly $\op$. 

Our next goal is to show that the groupoid of germs associated to the action of $\mathcal S$ on $\op$ is isomorphic to $\rp$ as a topological groupoid. In our case, $$Y = \{(s,T) \in \s \times \op :T\in D_{s^*s}\} =  \{([t_1,P,t_2],T) \in \s \times \op :T\in U(P,t_2)\},$$
and the equivalence relation on $Y$ can be characterized as follows:
\begin{lem}$\label{CharacEquivRel} ([t_1,P,t_2],T) \sim ([t_1',P',t_2'],T')$ if and only if $T = T'$, $t_1=t_1'$, $t_2=t_2'$ and there exists a patch $P''$ and a tile $t\in P''$, such that $P'' \supseteq P \cup P'$ and $T \in U(P'',t)$. 
\end{lem}

\begin{proof}
Let $s=[t_1,P,t_2]$, $s'=[t_1',P',t_2']$ and suppose $(s,T) \sim (s',T')$. Then $T=T'$ and there exists $s''=[t,P'',t]$ such that $T\in U(P'',t)$ and $s s''=s' s''$. Also, since $(s,T)$ and $(s',T')$ are in $Y$ we have that $T \in U(P,t_2)\cap U(P',t_2')\cap U(P'',t)$, which implies that $T(0)=t_2=t_2'=t$. It follows that $s s'' = [t_1, P \cup P'',t]$ and $s' s'' = [t_1', P' \cup P'',t]$ and hence $s s''=s' s''$ iff $P'' \supseteq P \cup P'$ and $t_1=t_1'$.

For the converse it is enough to take the idempotent as $e=[t, P'',t]$ in the definition of $\sim$.
\end{proof}
Note that if $(s,T)$ and $(s',T)$ are equivalent, then $x_s = x_{s'}$, and so both pairs represent translating $T$ by the same vector. This leads to the following.

\begin{theo}
The groupoid of germs associated to the action of $\mathcal S$ on $\op$ and $\rp$ are isomorphic as topological groupoids.
\end{theo}
\begin{proof}

Let $s=[t_1,P,t_2], s'=[t_1',P',t_2'] \in \s$. We will show that the map $$\alpha: \begin{array}{lll} \G & \rightarrow & \rp \\  $[s,T]$  & \mapsto & (T-x_s,T) \end{array},$$ where $x_s$ is the vector that takes the puncture of $t_2$ to the puncture of $t_1$, is a homeomorphism.

That $\alpha$ is well defined follows directly from lemma \ref{CharacEquivRel}. To see that $\alpha$ is a bijection note that its inverse is given by $\alpha^{-1}\left((T+x,T)\right) =[(T(0),P,T(x)),T],$ where $P$ is a patch containing $T(0) \cup T(x)$. 

Next we would like to show that the groupoid operations are preserved by the map $\alpha$. Recall that $\G^{\left(2 \right)}=\{ \left( [s,T],[s',T']\right) \in \G \times \G: T = \theta_{s'}(T')\}$, that is $\G^{\left(2 \right)}=\{ \left( [s,T],[s',T']\right) \in \G \times \G: T = T' - x_{s'}\}$ and notice that since $T\in U(P,t_2)$ and $T'\in U(P',t_2')$, $T = T' - x_{s'}$ implies that $t_2=t_1'$. Also recall that $\rp^{\left(2 \right)}=\{ \left( (T,T+x),(T',T'+x')\right) \in \rp \times \rp: T+x = T'\}$. Now it is not hard to check that $\left( [s,T],[s',T']\right) \in \G^{\left(2 \right)}$ if and only if $\left(\alpha([s,T]),\alpha([s',T'])\right)=\left((T-x_s,T),(T'-x_{s'},T') \right) \in \rp^{\left(2 \right)}$ and we have that:
$$\alpha([s,T][s',T']) = \alpha([ss',T'])= \alpha([[t_1,P\cup P', t_2'],T']) = (T',T' - x_r),$$ where $x_r$ the the vector from the puncture of $t_2'$ to the puncture of $t_1$. On the other hand 
$$\alpha([s,T])\alpha([s',T'])=(T-x_s,T)(T'-x_{s'},T')=(T'-x_s-x_{s'},T')=(T'-x_r,T')$$ and hence $\alpha$ is multiplicative.

To see that $\alpha$ preserves the inverse operation note that $$\alpha([s,T]^{-1})= \alpha([s^*,\theta_s(T)]) = \alpha([[t_2,P,t_1], T - x_s])=(T,T-x_s) = \alpha([s,T])^{-1}$$

For the topological part, notice that a basis for the topology of $\rp$ is given by the collection of graphs of the maps $s:U(P,t_2) \rightarrow U(P,t_1)$ given by $s(T)=T-x_s$ and a basis for the topology of $\G$ consists or the collection of all $\Theta(s,U)=\{[s,T]\in \G: T \in U\}$, where $s \in \s$ and $U$ is any open subset of $U(P,t_2)$. So it is clear that a basis element of $\rp$ is taken by $\alpha$ to an open set in $\G$. Now, let $\Theta(s,U)$ be a basis element of the the topology of $\G$. Then $\alpha(\Theta(s,U)) = \{ (T-x_s,T): T\in U \subseteq U(P,t_2) \}$. We show that $\alpha(\Theta(s,U))$ is open. Let $T_1 \in U$. Since $U$ is open there exists a patch $P'$ such that $T_1 \in U(P',t_2) \subseteq U \subseteq U(P,t_2)$ and this implies that $\{ (T-x_s,T): T \in U(P',t_2)\}$ is an open neighboorhood of $(T_1-x_s,T_1)$ (since it is the graph of $(t_2,P',t_1)$) contained in $\alpha(\Theta(s,U))$. It follows that $\alpha$ is bicontinuous and we conclude that $\alpha$ is a homeomorphism as desired.
\end{proof}

\section{Tight Representations of $\mathcal S$}
Recall that $E = \{ e\in \mathcal S \mid e^* = e^2 = e\}$, the set of idempotents in $\mathcal S$. In our case, they are elements of the form $[t, P, t]$. There is a natural partial order on $E$ given by $e\leq f \Leftrightarrow e = ef$. In our case, 
\[
[t, P, t] \leq [t', P', t'] \Leftrightarrow P'\subset P, t= t'
\]
for suitible choices of representatives of each class. We also notice that $[t, P, t] \leq [t', P', t']$ if and only if $t = t'$ and $U(P,t)\subset U(P', t')$.
\begin{defn}
Let $X$ be a partially ordered set with minimum element 0. A {\bf filter} is a non-empty subset $\xi\subset X$ such that 
\begin{enumerate}
\item $0\notin \xi$
\item if $x\in \xi$ and $y\geq x$, then $y\in\xi$.
\item if $x,y\in \xi$, $\exists z\in \xi$ such that $x, y \geq z$. 
\end{enumerate}
An {\bf ultra-filter} is a filter that is not properly contained in another. 
\end{defn}
We notice that filters in $E$ are closed under the product. Let $T$ be a partial tiling contained in some tiling in $\op$ and whose support contains the origin, and let $t = T(0)$. Then define
\[
\xi_T = \{ [t,P,t]\in E\mid t\in P\subset T\}.
\]
In words, $\xi_T$ is the set of all the $[t,P,t]$ such that $P$ is a patch around the origin in $T$. 
\begin{lem}\label{tilingultrafilter}
Let $T$ be a partial tiling contained in some tiling in $\op$ and whose support contains the origin. Then the set $\xi_T$ is a filter in $E$. If $T$ is a tiling, then $\xi_T$ is an ultra-filter. Furthermore, if $\xi\subset E$ is a non-empty filter, then there is a partial tiling $T$ such that $\xi = \xi_T$, and a tiling $T'\in\op$ such that $T\subset T'$. If $\xi$ is an ultra-filter, $T$ is neccessarily a tiling.
\end{lem}
\begin{proof} 

First we prove that $\xi_T$ is a filter: The zero element is clearly not in $\xi_T$. If we have $[t, P, t] \geq [t, P', t]$, and $[t, P', t]\in \xi_T$, then $P\subset P'$, and so $P\subset T$, hence $[t, P, t]\in \xi_T$. If $[t, P, t], [t, P', t]\in \xi_T$, then $P$ and $P'$ are compatible so $[t, P, t][t, P', t] = [t, P\cup P', t]$ is in $\xi_T$, and this product is less than both.

To see that if $T$ is a tiling then $\xi_T$ is an ultra-filter, suppose that we had a filter $\xi$ with $\xi_T \subsetneq \xi$. Take $[t', P', t']\in \xi \setminus \xi_T$, and find $M>0$ such that supp$(P) \subset B_M(0)$ (this is the ball of radius $M$ in $\rn$). Then $P = T(B_M(0))$ is a patch in $T$ around the origin, so $[t, P, t]\in \xi_T\subset \xi$. Filters are closed under products and 0 is not in a filter so $[t, P, t][t', P', t']$ must be non-zero. Hence $t = t'$ and $[t,P,t] \leq [t',P',t']$ implies that $[t',P',t']\in \xi_T$, a contradiction. We conclude that $\xi_T$ is an ultra-filter.

Now, suppose that $\xi$ is a non-empty filter. Then we can find $[t,P,t]\in \xi$. Filters in $E$ are closed under products, so if we have another element in $\xi$ it is necessarily of the form $[t, P', t]$ for some patch $P'$ that's compatible with $P$ and such that $U(P,t)\cap U(P',t) \neq \emptyset$. So let
\[
T = \bigcup_{[t,P',t]\in \xi}P'.
\]
As defined $T$ is a partial tiling because all such patches must be compatible since filters are closed under products and filters do not contain 0. It is clear from the definition of $\xi_T$ that $\xi\subset \xi_T$. Suppose that $[t, P', t]\in \xi_T$. For each tile $q\in P'$, there is a patch $P_q$ such that $q\in P_q$ and $[t, P_q, t]\in \xi$ because $T$ is the union of tiles in such patches. From this, we have that 
\[
P' \subset \bigcup_{q\in P'}P_q.
\]
Since $\xi$ is closed under products, $\left[t, \bigcup_{q\in P'}P_q, t\right]\in\xi$, and because of the above,
\[
[t, P', t] \geq \left[t, \bigcup_{q\in P'}P_q, t\right]
\]
and hence $[t, P', t] \in \xi$. It follows that $\xi=\xi_T$  as required. That $T$ is contained in some tiling $T'\in \op$ is obtained by noticing that the collection of closed sets $\{U(P',t)\}_{[t,P',t]\in \xi}$ has the property that any finite subcollection has non-empty intersection since all finite products are non-zero. Since $\op$ is compact, this implies that 
\[
\bigcap_{[t,P',t]\in \xi}U(P',t) \neq\emptyset.
\]
Now, if we take any $T'$ in the above intersection, then $T'$ agrees with $T$ everywhere, hence $T\subset T'$.

Now suppose that $\xi$ is an ultra-filter. From above, $\xi = \xi_T$ for some partial tiling $T$, and $T\subset T'$ for some tiling $T'$. It is clear that $\xi_T \subset \xi_{T'}$ and that the containment is proper unless $T$ is a tiling. 
\end{proof}

A {\em character} is a non-zero map
\[
\phi: E \to \{0,1\} \ \text{ such that } \ \phi(ef) = \phi(e)\phi(f)
\]
and which takes 0 to 0. Let $\hat{E}_0$ be the set of all characters on $E$. We give $\hat{E}_0$ the topology of pointwise convergence. If $\phi$ is a character, then the set
\[
\xi_\phi = \{ x\in E\mid \phi(x) = 1\}
\]
is a filter. Likewise, if $\xi$ is a filter, then the map
\[
\phi_\xi(x) =\left\{\begin{array}{ll}1,& \text{if } x\in \xi\\0,&\textnormal{otherwise}\end{array}\right.
\]
is a character. Characters on $E$ and filters in $E$ are in one-to-one correspondence. Let $\hat{E}_\infty$ be the set of all characters coming from ultra-filters, and denote its closure in $\hat{E}_0$ as $\hat{E}_{tight}$. Denote the character associated with $\xi_T$ as $\phi_T$. That is to say,
\[
\phi_T([t,P,t]) = \left\{\begin{array}{ll}1,& \text{if } T\in U(P,t)\\0, &\textnormal{otherwise}\end{array}\right. .
\]
Then we have the following
\begin{theo}
The map
\[
\Psi: \op \to \hat{E}_\infty
\]
\[
T \mapsto \phi_T
\]
is a homeomorphism. Furthermore, compactness of $\op$ means that $\op$ is homeomorphic to $\hat{E}_{tight}$ as well.
\end{theo} 
\begin{proof} Injectivity, at least intuitively, makes a lot of sense. The ultra-filter $\xi_T$ is the set of all patches around the origin in $T$. The union of all such patches is precisely $T$. We let $T\neq T'$ and find patches $P, P'\subset T'$ and $U\subset \rn$ around the origin such that $P$ and $P'$ do not overlap. 

Hence $[T(0), P, T(0)][T'(0),P', T'(0)] = 0$. Since filters are closed under the product, we must have that $\xi_T \neq \xi_{T'}$. Filters and characters are in one-to-one correspondence, so we must have that $\phi_T \neq \phi_{T'}$. Surjectivity is given by Lemma \ref{tilingultrafilter}.

Now, the topology on $\hat{E}_0$ has a sub-basis consisting of sets of the form
\[
b([t, P, t], U) = \{f\in \hat{E}_0\mid f([t, P, t])\in U\},
\] 
where $U$ is a subset of $\{0,1\}$. Consider the basic open set $U(P,t)\subset \op$. Then
\begin{eqnarray*}
\Psi(U(P,t)) &=& \{\phi_T \mid T\in U(P,t)\}\\
             &=& b([t, P, t], \{1\})\cap \hat{E}_\infty
\end{eqnarray*}
Hence $\Psi$ is an open map. Now 
$$\Psi^{-1}(b([t, P, t], \{0\})\cap \hat{E}_\infty) = U(P,t)^c$$
which is open due to $U(P,t)$ being closed. Hence $\Psi$ is bicontinuous and thus a homeomorphism. Since $\op$ is compact, $\hat{E}_\infty$ must be as well and so $\op\cong\hat{E}_\infty = \overline{\hat{E}_\infty} = \hat{E}_{tight}$. 
\end{proof}

Following \cite{\Exel} we know that there is an intrinsic action, called $\theta$, on any inverse semigroup $\s$ on $\hat{E}_{tight}$, such that for every idempotent $e\in E$, the domain of $\theta_e$ is $D_e=\{ \phi \in \hat{E}_{tight}: \phi(e)=1 \}$ and $\theta_s: D_{s^*s} \rightarrow D_{ss^*}$ is given by $\theta_s(\phi)(e)= \phi(s^*es)$. The groupoid of germs arising from this action is called the {\bf tight groupoid} of $\s$ and is denoted $\mathcal{G}_{\textnormal{tight}}$. Its associated C*-algebra, $C^*(\mathcal{G}_{\textnormal{tight}})$ is called the {\bf tight C*-algebra} of $\s$.

Recall that $\s$ also acts on $\op$ by the action $\theta^{\Omega}$, that is, for $s=[t_1,P,t_2]\in \s$, $\theta^{\Omega}_s:U(P,t_2) \rightarrow U(P,t_1)$ is defined via $\theta^{\Omega}_s(T) = T - x_s$, where $x_s$ is the vector from the puncture of $t_2$ to the puncture of $t_1$. Next we show that these two actions commute via the homeomorphism $\Psi$.

\begin{prop} For all $s=[t_1,P,t_2]\in \s$ the diagram below commutes. 
\begin{displaymath}
\xymatrix{
U(P,t_2) \ar[r]^{\theta^{\Omega}_s} \ar[d]_{\Psi} &
U(P,t_1) \ar[d]^{\Psi} \\
D_{s^*s} \ar[r]_{\theta_s} & D_{ss^*} }
\end{displaymath}

\end{prop}
\begin{proof}First we show that $\Psi(U(P,t_2))=\{\phi \in \hat{E}_{\infty}: \phi[t_2,P,t_2]=1\}$: Let $\phi \in \hat{E}_{\infty}$ be such that $\phi[t_2,P,t_2]=1$. By lemma \ref{tilingultrafilter}, $\phi = \phi_{\xi_T}$ for some $T \in \op$ and since $\phi_{\xi_T}[t_2,P,t_2] =1 $ iff $T \in U(P,t_2)$, we have that $\phi = \Psi (T)$. The other inclusion is trivial. Now, since $\psi$ is a homeomorphism, $\psi(U(P,t_2))$ is closed and hence equal to the closure of $\{\phi \in \hat{E}_{\infty}: \phi[t_2,P,t_2]=1\}$, which in turn is equal to $\{\phi \in \hat{E}_{tight}: \phi[t_2,P,t_2]=1\}$, since if $\phi \in \hat{E}_{tight}$ and $\phi[t_2,P,t_2]=1$ then there exists a sequence $(\phi_i)$ in $\hat{E}_{\infty}$ converging to $\phi$, which implies that for $i$ large enough $\phi_i[t_2,P,t_2]=\phi[t_2,P,t_2]=1$. It follows that $\Psi(U(P,t_2))= \{\phi \in \hat{E}_{tight}: \phi[t_2,P,t_2]=1\} = D_{s^*s}$. Analogously one can show that $\Psi(U(P,t_1))=D_{ss^*}$.

Now let $T\in U(P,t_2)$. Then $$\theta_s(\Psi(T))[t,P',t]=\theta_s(\phi_t)[t,P',t]= \phi_T \left( [t_2,P,t_1][t,P',t][t_1,P,t_2]\right)$$
and the element in which $\phi_T$ is applied is non-zero iff $t=t_1$ and $U(P,t_1)\cap U(P',t_1) \neq \emptyset$, which happens if and only if $t=t_1$ and $P$ and $P'$ are compatible, in which case $$\theta_s(\Psi(T))[t,P',t]=\phi_T \left( [t_2,P \cup P',t_2] \right).$$ It follows that $$\theta_s(\Psi(T))[t,P',t]=  \begin{cases}1,& \text{if } t=t_1 \text{ and $P$ and $P'$ are compatible.} \\0, &\textnormal{otherwise}\end{cases}.$$ 

On the other hand, $$\Psi( \theta^{\Omega}_s(T)) [t,P',t]=\Psi( T-x_s) [t,P',t]= \phi_{T-x_s}[t,P',t]= \begin{cases} 1,& \text{if } T-x_s \in U(P',t) \\0, &\textnormal{otherwise}\end{cases}$$
and since $T\in U(P,t_2)$ we have that $T-x_s \in U(P',t)$ iff $t=t_1$ and $P$ and $P'$ are compatible. It follows that $\theta_s\circ \Psi = \Psi \circ \theta^{\Omega}_s$ as desired.
\end{proof}

We combine these to obtain our main result.

\begin{theo}
Let $T$ be a tiling which is strongly aperiodic, repetitive, and has finite local complexity and let $\s$ be its associated inverse semigroup. Then the groupoids $\rp$ and tight groupoid associated with $\mathcal S$, $\mathcal G_{\textnormal{tight}}$ are isomorphic as topological groupoids. In particular, the tight C*-algebra associated with $\mathcal S$ is isomorphic to the usual C*-algebra associated to a tiling $C^*(\rp)$.
\end{theo}

\end{document}